\documentclass[12pt]{article}
\usepackage{amsfonts}
\usepackage{latexsym}
\usepackage{cite}
\usepackage{amsmath,amsfonts,latexsym,amssymb,color}
\usepackage[mathscr]{eucal}
\usepackage{cases}
\usepackage{amsthm}

\usepackage[bf,small]{caption2}
\usepackage{float}
\usepackage{graphicx}
\usepackage{amsmath}
\usepackage{amssymb}
\usepackage[all]{xy}

\newtheorem{theorem}{theorem}[section]
\newtheorem{thm}[theorem]{Theorem}
\newtheorem{lem}[theorem]{Lemma}

\newtheorem{prop}[theorem]{Proposition}

\newtheorem{rmk}[theorem]{Remark}
\newtheorem{nota}[theorem]{Notation}

\begin{document}

\title{\textbf{Lifting graph automorphisms along solvable regular covers}}
\author{\Large Haimiao Chen
\footnote{Email: \emph{chenhm@math.pku.edu.cn}}\\
\normalsize \em{Mathematics, Beijing Technology and Business
University, Beijing,
100037, P.R. China}\\
\Large Jin Ho Kwak
\footnote{Email: \emph{jinkwak@postech.ac.kr}}\\
\normalsize \em{Mathematics,  Beijing Jiaotong University, Beijing,
100044, P.R. China }\\
\normalsize \em{Mathematics, POSTECH, Pohang, 790-784, Republic of
Korea}}

\date{}
\maketitle

\begin{abstract}
  A {\em solvable} cover of a graph is a regular cover whose covering transformation group is solvable. In this paper, we show that a solvable cover of a graph can be decomposed into layers of abelian covers, and also, a lift of a given automorphism of the base graph of a solvable cover can be decomposed into layers of lifts of the automorphism in the layers of the abelian covers.
  This procedure is applied to classify metacyclic covers of the tetrahedron branched at face-centers.

\end{abstract}


\section{Introduction}

In this paper all graphs are finite and connected.
For a graph $\Gamma$, let $V(\Gamma)$, $E(\Gamma)$ and $D(\Gamma)$ denote the set of vertices, edges and arcs of $~\Gamma$, respectively.
For a vertex $v\in V(\Gamma)$, let $N(v)$ be the set of arcs emanating from $v$, called the {\em neighborhood} of $v$.
For convenience, we consider each vertex $v$ as a {\em trivial arc} at $v$, which means an arc of length 0 from $v$ to itself. Let $\check{D}(\Gamma)=V(\Gamma)\sqcup D(\Gamma)$, and let $\iota: V(\Gamma)\hookrightarrow\check{D}(\Gamma)$ denote the inclusion.
For each $x\in \check{D}(\Gamma)$, let $s(x)$ and $t(x)$ denote its starting vertex and terminal vertex, respectively; this defines two functions $s,t:\check{D}(\Gamma)\to V(\Gamma)$. Call the triple $(\iota,s,t)$ the \emph{incidence relation} of $~\Gamma$.
Each edge $e\in E(\Gamma)$ gives rise to two arcs $x, x'$ which are opposite to each other; we say that $e$ is the \emph{underlying edge} of $x$ and $x'$.
A \emph{walk} $W=(x_{1},\ldots,x_{n})$ is a sequence of arcs such that $t(x_{i})=s(x_{i+1})$ for $1\leqslant i<n$; in addition, if $s(x_{1})=t(x_{n})=v$, we call
$W$ a \emph{closed walk} at $v$.
Let $T$ be a spanning tree of a graph $\Gamma$  and let $u_0\in V(\Gamma)$ be a fixed base vertex. For each $u\in V(\Gamma)$, let $W(u)$ denote the unique reduced walk from $u_0$ to $u$ in the tree $T$. Choose an arc for each cotree edge, and let $x_{i},i=1,\ldots,\beta$ denote these cotree arcs, where $\beta=|E(\Gamma)|-|V(\Gamma)|+1$ is the Betti number of $~\Gamma$. For each $i$, let $L_{i}$ be the closed walk $W(s(x_{i}))\cdot x_{i}\cdot W(t(x_{i}))^{-1}$,  called a {\em fundamental closed walk}. Then  the {\em fundamental group} $\pi_{1}(\Gamma,u_{0})$ is the free group generated by $L_{1},\ldots,L_{\beta}$.

A \emph{graph morphism} $f:\Gamma\to\Gamma'$ is a function $f:\check{D}(\Gamma)\to\check{D}(\Gamma')$ preserving the incidence relation.
A \emph{cover} $\pi:\tilde{\Gamma}\rightarrow\Gamma$ is a graph morphism such that $\pi:\check{D}(\tilde{\Gamma})\rightarrow \check{D}(\Gamma)$ is surjective
and $\pi|_{N(\tilde{v})}:N(\tilde{v})\rightarrow {N(v)}$
is bijective for any $v\in V(\Gamma)$ and $\tilde{v}\in \pi^{-1}(v)$. Call the graph $\tilde{\Gamma}$ the {\it covering graph}, call $\Gamma$ the
{\it base graph}, and call $\pi^{-1}(v)$ the {\it fiber over $v$}. The group of automorphisms of the covering graph $\tilde{\Gamma}$ which fix each fiber setwise is called the \emph{covering transformation group} of $\pi$ and denoted by ${\rm CTG}(\pi)$. The cover $\pi$ is called \emph{regular} if ${\rm CTG}(\pi)$ acts transitively on each fiber, and usually called an \emph{$A$-cover} when ${\rm CTG}(\pi)$ is isomorphic to a finite group $A$.

For two covers $\pi:\tilde{\Gamma}\rightarrow\Gamma$ and $\pi':\tilde{\Gamma}'\rightarrow\Gamma'$, an \emph{isomorphism} $\pi\rightarrow\pi'$
is a pair of isomorphisms $(\alpha,\tilde{\alpha})$ consisting of $\alpha:\Gamma\to\Gamma'$ and $\tilde{\alpha}:\tilde{\Gamma}\rightarrow\tilde{\Gamma}'$
satisfying $\alpha\circ\pi=\pi'\circ\tilde{\alpha}$. Call $\tilde{\alpha}$ a \emph{lift} of $\alpha$ along $(\pi,\pi')$, and say that
$\alpha$ is lifted to $\tilde{\Gamma}\rightarrow\tilde{\Gamma}'$. When $\Gamma=\Gamma'$ and $\alpha$ is the identity morphism, call $\tilde{\alpha}$ an \emph{equivalence} between covers and say that $\pi$ is \emph{equivalent} to $\pi'$.
When $\pi=\pi'$, simply call $\tilde{\alpha}$ a lift of $\alpha$ along $\pi$, and say that $\alpha$ is lifted to $\tilde{\Gamma}$.

Given a finite group $A$, a \emph{voltage assignment} of $~\Gamma$ in $A$ is a function $\phi:D(\Gamma)\rightarrow A$ such that
$\phi(x^{-1})=\phi(x)^{-1}$ for all $x\in D(\Gamma)$.
Each voltage assignment $\phi$ determines an $A-$cover $\pi_{\phi}:\Gamma\times_{\phi}A\to\Gamma$, called the {\it derived cover} of $\phi$, as follows. The covering graph $\Gamma\times_{\phi}A$ has vertex set $V(\Gamma\times_{\phi}A)=V(\Gamma)\times A$ and arc set $D(\Gamma\times_{\phi}A)=D(\Gamma)\times A$ with $s\langle x,g\rangle=(s(x),g), t\langle x,g\rangle=(t(x),g\phi(x))$; here we adopt the notational convention that vertices (trivial arcs) of $\Gamma\times_{\phi}A$ are denoted by pairs $(u,g)$ while arcs are denoted by $\langle x,g\rangle$. In fact, $\pi_{\phi}:\Gamma\times_{\phi}A\to\Gamma$ is the canonical morphism given by the projection onto the first coordinate, and the action of ${\rm CTG}(\pi)$ is given by the canonical left $A$-action
$(h,\langle x,g\rangle)\mapsto\langle x,hg\rangle=:h\langle x,g\rangle.$

Let $T$ be a spanning tree. A voltage assignment $\phi:D(\Gamma)\rightarrow A$ is said to be {\it $T$-reduced} if $\phi(x)=1$ for all $x\in D(T)$.
It is a classical result (see \cite{GT77}) that each $A$-cover of $\Gamma$ is equivalent to the derived cover of some voltage assignment; furthermore, this voltage assignment can be chosen to be $T$-reduced for any previously given spanning tree $T$.

For any walk $W=(x_{1},\ldots,x_{n})$, the voltage of $W$ is naturally defined by $\phi(W)=\phi(x_{1})\cdots\phi(x_{n})$.

Lifting problems are traditionally studied within the context of regular
covers, specially by voltage assignments; see \cite{Ma98, MNS00, MMP-cover2004}. Here we recall some key ingredients.

Let $\tilde{\Gamma}_{i}=\Gamma_{i}\times_{\phi_{i}}A, i=1,2$ be two covers and let an isomorphism $\alpha:\Gamma_{1}\to\Gamma_{2}$  be lifted to $\tilde{\alpha}:\tilde{\Gamma}_{1}\to\tilde{\Gamma}_{2}$ which is a function $\tilde{\alpha}:\check{D}(\Gamma_{1})\times A\to \check{D}(\Gamma_{2})\times A$.
Let $\xi:\check{D}(\Gamma_{1})\times A\to A$ denote the composition of $\tilde{\alpha}$ with the projection $\check{D}(\Gamma_{2})\times A\to A$ so that
\begin{align}\label{eq=def(xi)}
\tilde{\alpha}\langle x,g\rangle=\langle\alpha(x),\xi\langle x,g\rangle\rangle, \hspace{5mm} x\in \check{D}(\Gamma_{1}), g\in A.
\end{align}
Considering its starting and terminal vertices, we get
\begin{align*}
\tilde{\alpha}(s(x),g)=\tilde{\alpha}(s\langle x,g\rangle)=s(\tilde{\alpha}\langle x,g\rangle)=s\langle\alpha(x),\xi\langle x,g\rangle\rangle=(s(\alpha(x)),\xi\langle x,g\rangle), \\
\tilde{\alpha}( t(x),g\phi_{1}(x))=t(\tilde{\alpha}\langle x,g\rangle)=t\langle\alpha(x),\xi\langle x,g\rangle\rangle=(t(\alpha(x)),\xi\langle x,g\rangle\phi_{2}(\alpha(x))).
\end{align*}
The first equation shows that $\xi\langle x,g\rangle$ depends only on $s(x)$ and $g$, but not on the terminal vertex $t(x)$;
in particular,
\begin{align}
\xi\langle x,g\rangle=\xi(s(x),g). \label{eq:xi2}
\end{align}
Now, by using the second equation and the definition of $\xi$, one can see that
$$\xi(t(x),g\phi_{1}(x))=\xi(s(x),g)\phi_{2}(\alpha(x)),$$
which implies
\begin{align}\label{eq=xi W}
\xi(v,g\phi_{1}(W))=\xi(u,g)\phi_{2}(\alpha(W))
\end{align}
for any walk $W$ from $u$ to $v$ and for any $g\in A$.

Choose a vertex $u_0$ as a base in the graph $\Gamma_1.$ For any $(u,g)$, take a walk $W$ from $u_0$ to $u$ with $\phi_{1}(W)=g$, (this is always possible because for any value $g'\in A$, the connectedness assumption assures the existence of a closed walk $W'$ at $u_0$ with $\phi_{1}(W')=g'$). Then, by Eq.~(\ref{eq=xi W})
\begin{align}
\xi(u,g)=\xi(u,1\phi_{1}(W))=\xi(u_0,1)\phi_{2}(\alpha(W)), \label{eq:xi}
\end{align}
where $1$ is the identity in the voltage group $A$.
In particular, $\xi$ is completely determined by $\alpha$ and $\xi(u_{0},1)$.

Conversely, if a function $\xi:\check{D}(\Gamma_{1})\times A\to A$ is defined by Eqs.~(\ref{eq:xi}) and (\ref{eq:xi2}) with a given value $\xi(u_0,1)$ in $A$,
then it gives rise to a lift $\tilde{\alpha}$ of $\alpha$ by Eq.~(\ref{eq=def(xi)}). A sufficient and necessary condition for such $\xi$ (so also $\tilde{\alpha}$) to be well-defined is that, for any closed walk $W$ at $u_0$, whenever $\phi_{1}(W)=1$, we have $\phi_{2}(\alpha(W))=1$.
This can be rephrased as the following proposition which is a generalization of \cite{HKL96} Theorem 3:

\begin{prop} \label{prop:lifting}
An isomorphism $\alpha:\Gamma_{1}\to\Gamma_{2}$ can be lifted to $\Gamma_{1}\times_{\phi_{1}}A\to\Gamma_{2}\times_{\phi_{2}}A$ if and only if there
exists an automorphism $\sigma:A\to A$ such that the following diagram commutes:
$$
\xymatrix{
\pi_{1}(\Gamma_{1},u_{0})\ar[d]_{\phi_{1}}\ar[r]^{\alpha_{\ast}}&\pi_{1}(\Gamma_{2},\alpha(u_{0}))\ar[d]^{\phi_{2}} \\
A\ar[r]^{\sigma}&A
}
$$ 
Furthermore, for any two lifts $\tilde{\alpha}$ and $\tilde{\alpha}'$, there exists a unique (covering transformation) $h\in A$ such that
$$\tilde{\alpha}'\langle x,g\rangle=h\tilde{\alpha}\langle x,g\rangle \text{\ for\ all\ } x\in \check{D}(\Gamma_{1}), \ g\in A.$$
\end{prop}

\bigskip

The lifting problem for graph covers is of fundamental importance and has attracted much attention recently; see \cite{HKL96, Ma98, MNS00, Sk86}.
Two regular covers of a graph $\Gamma$ are equivalent if and only if the identity morphism of $\Gamma$ can be lifted, and they are isomorphic if and only if at least one automorphism of $~\Gamma$ can be lifted. Usually, to construct large graphs with high symmetries, people like to construct {\it admissible covers}, which mean regular covers that admit lifts of certain given automorphisms of the base graph.
Moreover, some problems in map theory are closely related to admissible covers of graphs.

Du et al. \cite{DKX03} gave a linear criterion for lifting automorphisms along elementary abelian covers, and applied it to classify arc-transitive elementary  abelian covers of the Petersen graph. In the work of Malni\v c and Poto\v cnik \cite{MP06}, a method of ``invariant subspace" was developed
to deal with elementary abelian covers. For general abelian covers, the first author \cite{Ch13} proposed a useful criterion, as
summarized in the next section.

Besides the above, till now, there is no much work on the lifting problem for nonabelian covers. Malni\v c et al.
\cite{MMP2004}, which dealt with cubic graphs admitting a solvable edge transitive group, proved that such graphs can be obtained by successive
elementary abelian covers of the 3-dipole $D_{3}$ or the complete graph $K_{4}$. Recently Xu et al. \cite{XDKX14} classified 2-arc-transitive metacyclic covers of complete graphs.
In this paper we  discuss the nonabelian case when the voltage group $A$ is presented as an extension of a characteristic subgroup $N$ by a group $Q$,
then each $A$-cover $\Gamma\times_{\phi}A\to\Gamma$ can be decomposed into a composition of a $Q$-cover $\Gamma\times_{\overline{\phi}}Q\to\Gamma$
and an $N$-cover $\Gamma\times_{\phi}A\to\Gamma\times_{\overline{\phi}}Q$
with an induced voltage assignment $\overline{\phi}$ from ${\phi}$. Moreover, an automorphism $\alpha\in\textrm{Aut}(\Gamma)$ can be lifted to $\Gamma\times_{\phi}A$ if and only if it can be lifted to $\Gamma\times_{\overline{\phi}}Q$ and all of its lifts are further lifted to $\Gamma\times_{\phi}A$. Consequently all (connected) covers with solvable covering transformation groups can be decomposed into layers of abelian
covers, and an automorphism can be lifted if and only if it can be lifted  step by step through the layers upwards.

This paper is organized as follows. In Section 2 we quickly review the criterion given by \cite{Ch13}, for lifting automorphisms along abelian covers.
In Section 3 we propose a criterion for lifting automorphisms along a kind of nonabelian covers, and then give a criterion for solvable covers. In Section 4, as an application, we classify metacyclic covers of the tetrahedron.

As a notational convention, we use $\mathbb{Z}_{n}$ to denote the cyclic group $\mathbb{Z}/n\mathbb{Z}$, and also consider it as a quotient ring of $\mathbb{Z}$ whenever it is necessary.

\section{Review of the theory of lifts along abelian covers}

Let $\Gamma$ be a graph with a spanning tree and a base vertex $u_0\in V(\Gamma)$ chosen, and let $L_{1},\ldots,L_{\beta}$ be the fundamental closed
walks in  $\Gamma$. The {\em homology group} $H_{1}(\Gamma;\mathbb{Z})$  of $~\Gamma$ is defined as the abelianized group of the  fundamental group $\pi_{1}(\Gamma,u_{0})$, and hence it is the free abelian group generated by the  fundamental closed walks $L_{1},\ldots,L_{\beta}$, implying
$H_{1}(\Gamma;\mathbb{Z})\cong\mathbb{Z}^{\beta}.$

Let $A=\prod_{\gamma=1}^{g}A_{\gamma}$ be any finite abelian group, with distinct prime divisors $p_{\gamma}, \gamma=1,\ldots,g$  and the $p_{\gamma}$
subgroups
\begin{align}
A_{\gamma}=\prod
\limits_{\eta=1}^{n_{\gamma}}\mathbb{Z}_{p_{\gamma}^{k(\gamma,\eta)}}, \hspace{5mm}
k(\gamma,1)\leqslant\cdots\leqslant k(\gamma,n_{\gamma}).
\end{align}

Suppose $\phi:D(\Gamma)\to A$ is a voltage assignment. It extends to a group homomorphism which we will also denote by
\begin{align}
\phi:H_{1}(\Gamma;\mathbb{Z})\to A,
\end{align}
because $A$ is assumed to be abelian.
In order to restrict ourselves to the case of connected A-covers, we will always assume that it is surjective.

For $\gamma\in\{1,\ldots, g\}$, let $\sigma_{\gamma}:A\twoheadrightarrow A_{\gamma}$ and $\tau_{\gamma}:A_{\gamma}\hookrightarrow A$ denote the canonical
projection and injection in groups, respectively, and let $R_{\gamma}=\mathbb{Z}_{p_{\gamma}^{k(\gamma, n_{\gamma})}}$.
For each $\eta\in\{1,\ldots, n_{\gamma}\}$, there is an injection (as a ring homomorphism)
\begin{align}
\mathbb{Z}_{p_{\gamma}^{k(\gamma,\eta)}}\hookrightarrow R_{\gamma}, \hspace{5mm}
\lambda\mapsto p_{\gamma}^{k(\gamma,n_{\gamma})-k(\gamma,\eta)}\cdot\lambda.
\end{align}

The $p_{\gamma}$ parts $\sigma_{\gamma}(\phi(L_{i}))$ of $\phi(L_{i})$ shall be interpreted as row vectors of a module ${R_{\gamma}}^{n_{\gamma}}$ over the ring $R_{\gamma}$, no longer in the group $A_{\gamma}$, and placed in rows in order to form a nonzero $\beta\times n_{\gamma}$ matrix $[\phi]_{\gamma}$
with entries in $R_{\gamma}$.

A graph automorphism $\alpha\in\textrm{Aut}(\Gamma)$ induces a group automorphism
\begin{align}
\alpha_{\ast}:H_{1}(\Gamma;\mathbb{Z})\rightarrow H_{1}(\Gamma;\mathbb{Z}),
\end{align}
represented by a matrix $[\alpha]=[\alpha_{i,j}]\in\textrm{GL}(\beta,\mathbb{Z})$ with respect to a basis $\{L_{1},\ldots,L_{\beta}$\}; it is given by
\begin{align}
\alpha_{\ast}(L_{i})=\sum\limits_{j=1}^{\beta}\alpha_{i,j}L_{j}.   \label{eq:matrix}
\end{align}

It is known (see \cite{Ch13, DKX03} for example) that the automorphism $\alpha$ can be lifted if and only if the null space $\mathcal{N}([\phi]_{\gamma})$
of the matrix $[\phi]_{\gamma}$ is invariant under the automorphism $\alpha_{\ast}$ for each ${\gamma}$. Moreover, the veracity of the sufficiency condition can be easily determined with
basis vectors  if a basis for the null space $\mathcal{N}([\phi]_{\gamma})$ is given.
To find a basis for the null space $\mathcal{N}([\phi]_{\gamma})$, we choose change of basis matrices $B_{\gamma}\in\textrm{GL}(\beta,R_{\gamma})$
and $C_{\gamma}\in\textrm{GL}(n_{\gamma},R_{\gamma})$ so that $B_{\gamma}[\phi]_{\gamma}C_{\gamma}$ is in normal form, i.e., $(B_{\gamma}[\phi]_{\gamma}C_{\gamma})_{i,j}=\delta_{i,j}\cdot p_{\gamma}^{s_{\gamma,i}}$ for all $i,j.$ Note that the matrix
$B_{\gamma}[\alpha]B_{\gamma}^{-1}$ is a representation of $\alpha_{\ast}$ with respect to a new basis by the change of basis matrix $B_{\gamma}$, where $[\alpha]$ should be understood as a matrix in $\textrm{GL}(\beta,R_{\gamma})$ through the quotient homomorphism
$\mathbb{Z}\twoheadrightarrow R_{\gamma}$.

The next theorem was shown by the first author \cite{Ch13}:
\begin{thm} \label{thm:abelian}
Under the notation presented above, $\alpha$ can be lifted if and only if $\deg_{p_{\gamma}}((B_{\gamma}[\alpha]B_{\gamma}^{-1})_{i,j})\geqslant s_{\gamma,i}-s_{\gamma,j}$ for all $\gamma$ and all $(i,j)$ with $1\leqslant j<i\leqslant \beta$.
Here, for an element $\mu\in R_{\gamma}$, $\deg_{p_{\gamma}}(\mu)$ is defined to be  the largest integer $r$ such that $p_{\gamma}^{r}$ divides $\mu$.
\end{thm}

\section{Lifting along nonabelian covers}

\subsection{Decomposing a regular $A$-cover}

Let $A$ be an extension of a subgroup
$N$ by a group $Q$, so there is a (split or not) short exact sequence of finite groups
$$1\to N\stackrel{\imath}{\rightarrow} A\stackrel{\jmath}{\rightarrow} Q\to 1.$$

Let $\phi:D(\Gamma)\rightarrow A$ be a voltage assignment.
Then the composite
\begin{align}
\overline{\phi}:=\jmath\circ\phi:D(\Gamma)\rightarrow Q  \label{eq:phibar}
\end{align}
is a voltage assignment in $Q$, 
and the function
\begin{align}\label{Eq=cover1}
\Gamma\times_{\phi}A\to\Gamma\times_{\overline{\phi}}Q, \hspace{5mm} \langle x,g\rangle\mapsto \langle x,\jmath(g)\rangle
\end{align}
defines an $N$-cover of $\Gamma\times_{\overline{\phi}}Q$.

Fix a function $\kappa:Q\to A$ such that $\jmath\circ\kappa$ is the identity on $Q$ and $\kappa(1)=1$. (Such a function $\kappa$ is called a \emph{section}; it is not necessarily a group homomorphism.)
By noting that the images of $\kappa(q)\phi(x)$ and $\kappa(q\overline{\phi}(x))=\kappa\jmath[\kappa(q)\phi(x)]$ under $\jmath$ are the same in $Q$,
one can define a voltage assignment
\begin{align}
\underline{\phi}:D(\Gamma\times_{\overline{\phi}}Q)\rightarrow N, \ \ \ \ \
\langle x,q\rangle\mapsto [\kappa(q)\phi(x)][\kappa(q\overline{\phi}(x))]^{-1}, \label{Eq=underlinephi}
\end{align}
so as to get a new cover of $\Gamma\times_{\overline{\phi}}Q$
\begin{align}\label{Eq=cover2}
(\Gamma\times_{\overline{\phi}}Q)\times_{\underline{\phi}}N \to \Gamma\times_{\overline{\phi}}Q.
\end{align}

Now the two $N$-covers $(\ref{Eq=cover1})$ and $(\ref{Eq=cover2})$ of $\Gamma\times_{\overline{\phi}}Q$ are equivalent.
In fact, an equivalence between them can be given by
\begin{align}
\Xi_{\phi}:\Gamma\times_{\phi}A\rightarrow(\Gamma\times_{\overline{\phi}}Q)\times_{\underline{\phi}}N,\hspace{5mm}
\langle x,g\rangle\mapsto \langle\langle x,\jmath(g)\rangle,g[\kappa(\jmath(g))]^{-1}\rangle;
\end{align}
this function indeed preserves the incidence relation because
$$g[\kappa(\jmath(g))]^{-1}\underline{\phi}\langle x,\jmath(g)\rangle=g\phi(x)[\kappa(\jmath(g\phi(x)))]^{-1}.$$
By taking compositions of $N$-cover projections $(\ref{Eq=cover1})$ and $(\ref{Eq=cover2})$ with the projection
$\Gamma\times_{\overline{\phi}}Q \to \Gamma$, one can have two equivalent $A$-covers of $~\Gamma$. Actually, the $A$-action on
$(\Gamma\times_{\overline{\phi}}Q)\times_{\underline{\phi}}N$ can be given by ``pushing forward" the action of $A$ on $\Gamma\times_{\phi}A$  via the equivalence $\Xi_{\phi}$, which explicitly means
$$(g,\langle\langle x,q\rangle,h\rangle)\mapsto \langle\langle x,\jmath(g)q\rangle,gh\kappa(q)[\kappa(\jmath(g)q)]^{-1}\rangle$$
for any $g\in A,  \langle x,q\rangle\in D(\Gamma\times_{\overline{\phi}}Q)$ and $ h\in N.$

In this sense, we have decomposed the original cover $\Gamma\times_{\phi}A\to \Gamma$ into the $N$-cover $(\Gamma\times_{\overline{\phi}}Q)\times_{\underline{\phi}}N\to\Gamma\times_{\overline{\phi}}Q$ and the $Q$-cover $\Gamma\times_{\overline{\phi}}Q\to\Gamma$.

Conversely, given two voltage assignments $\varphi: D(\Gamma)\to Q$ and $\psi:D(\Gamma\times_{\varphi}Q)\to N$, does there exist a voltage assignment $\phi: D(\Gamma)\to A$ such that $\overline{\phi}=\varphi$ and $\underline{\phi}=\psi$?
If these two conditions hold, then by (\ref{Eq=underlinephi}), the three voltage assignments $\varphi,\psi$ and $\phi$ satisfy
\begin{align}
\psi\langle x,q\rangle=\kappa(q)\phi(x)[\kappa(q\varphi(x))]^{-1}, \ \ \ \text{for\ all\ }x\in D(\Gamma), q\in Q. \label{Eq=relation}
\end{align}
In particular,
$$\psi\langle x,1\rangle=\kappa(1)\phi(x)[\kappa(1\varphi(x))]^{-1}=\phi(x)[\kappa(\varphi(x))]^{-1},$$
hence $\phi$ should be defined by
\begin{align} \label{eq=def(phi)}
\phi(x)=\psi\langle x,1\rangle\kappa(\varphi(x)).
\end{align}
Then (\ref{Eq=relation}) is equivalent to
\begin{align}\label{Eq=well psi}
\psi\langle x,q\rangle=\kappa(q)\psi\langle x,1\rangle\kappa(\varphi(x))[\kappa(q\varphi(x))]^{-1}.
\end{align}

Certainly, there exist some pairs of $\varphi$ and $\psi$ which do not satisfy (\ref{Eq=well psi}) and hence, to serve our purpose, we need to impose additional restrictions.
For a given  section $\kappa$ and for each $\varphi: D(\Gamma)\to Q$, let $\mathcal{C}_{\kappa}(\varphi;N)$ denote the set of voltage assignments $\psi:D(\Gamma\times_{\varphi}Q)\to N$ which satisfy (\ref{Eq=well psi}).

We thus have established the following two theorems.

\begin{thm}\label{Thm=cover-decompose}
Let $\phi:D(\Gamma)\to A$ be a voltage assignment, and let
$A$ be an extension of a subgroup $N$ by a group $Q$. Then $\phi$  induces two voltage assignments $\overline{\phi}:D(\Gamma)\to Q$ and $\underline{\phi}:D(\Gamma\times_{\overline{\phi}}Q)\to N$ such that the two covers
$$\xymatrix{
  \Gamma\times_{\phi}A \ar[rr]^{\Xi_{\phi}} \ar[dr]
                &  &    (\Gamma\times_{\overline{\phi}}Q)\times_{\underline{\phi}}N\ar[dl]    \\
                & \Gamma                 }
$$
are equivalent as $A$-covers.
\end{thm}

\begin{thm} \label{thm:bijection}
Given a fixed section $\kappa:Q\to A$ with $\kappa(1)=1$, there is a one-to-one correspondence between the set of voltage assignments
$\phi:D(\Gamma)\rightarrow A$ and the set of pairs of voltage assignments $(\varphi:D(\Gamma)\rightarrow Q,\psi\in\mathcal{C}_{\kappa}(\varphi;N))$.
\end{thm}

\begin{rmk} \label{rmk:T-reduced}
\rm From (\ref{eq:phibar}), (\ref{Eq=underlinephi}) and (\ref{eq=def(phi)}) we see that for each $x\in D(\Gamma)$, $\phi(x)=1$ if and only if $\varphi(x)=1$ and $\psi\langle x,1\rangle=1$. Thus, given a spanning tree $T$, to obtain a $T$-reduced voltage assignment $\phi:D(\Gamma)\to A$, it is equivalent to find a $T$-reduced $\varphi:D(\Gamma)\to Q$ and $\psi\in\mathcal{C}_{\kappa}(\varphi;N)$ with $\psi\langle x,1\rangle=1$ for all $x\in D(T)$.
\end{rmk}

The next theorem shows that a lifting problem on an $A$-cover can be decomposed into two lifting problems by using Theorem~\ref{Thm=cover-decompose}.

\begin{thm} \label{Thm=lift-decompose}
Let $A$ be an extension of a characteristic subgroup $N$ by a group $Q$ and let
$\phi_{i}: D(\Gamma_{i})\to A, i=1,2$ be voltage assignments. Then an isomorphism $\alpha:\Gamma_{1}\to\Gamma_{2}$ can be lifted to $\Gamma_{1}\times_{\phi_{1}}A\to\Gamma_{2}\times_{\phi_{2}}A$ if and only if
\begin{enumerate}
  \item[{\rm (i)}]  $\alpha$ can be lifted to $\Gamma_{1}\times_{\overline{\phi_{1}}}Q\to\Gamma_{2}\times_{\overline{\phi_{2}}}Q$, and
  \item[{\rm (ii)}]  any such lift can be lifted again to
$(\Gamma_{1}\times_{\overline{\phi_{1}}}Q)\times_{\underline{\phi_{1}}}N\to(\Gamma_{1}\times_{\overline{\phi_{2}}}Q)\times_{\underline{\phi_{2}}}N$.
\end{enumerate}
Furthermore, in the condition {\rm(ii)}, if one of such lifts can be lifted to
$\Gamma_{1}\times_{\overline{\phi_{1}}}Q)\times_{\underline{\phi_{1}}}N\to(\Gamma_{1}\times_{\overline{\phi_{2}}}Q)\times_{\underline{\phi_{2}}}N$,
then all other lifts can also be lifted.
\end{thm}

\begin{proof}
($\Rightarrow$) Let $\alpha$ be lifted to $\tilde{\alpha}: \Gamma_{1}\times_{\phi_{1}}A\to\Gamma_{2}\times_{\phi_{2}}A$. Then by Proposition
\ref{prop:lifting}, there exists an automorphism $\sigma:A\to A$ such that $\sigma\circ\phi_{1}=\phi_{2}\circ\alpha_{\ast}$. Since $N$ is assumed to be characteristic, $\sigma$ induces an automorphism $\overline{\sigma}:Q\to Q$ such that $\overline{\sigma}\circ\jmath=\jmath\circ\sigma$. Thus $\overline{\sigma}\circ\overline{\phi_{1}}=\overline{\phi_{2}}\circ\alpha_{\ast}$, and by Proposition \ref{prop:lifting} again, $\alpha$ can be lifted
to $\Gamma_{1}\times_{\overline{\phi_{1}}}Q\to\Gamma_{2}\times_{\overline{\phi_{2}}}Q$.

To show (ii), let $\overline{\alpha}:\Gamma_{1}\times_{\overline{\phi_{1}}}Q\to\Gamma_{2}\times_{\overline{\phi_{2}}}Q$ be any lift of $\alpha$ and let $\overline{W}$ be any closed walk in $\Gamma_{1}\times_{\overline{\phi_{1}}}Q$ with $\underline{\phi_{1}}(\overline{W})=1$. One can write
$$\overline{W}=(\langle x_{1},q\rangle,\langle x_{2},q\overline{\phi_{1}}(x_{1})\rangle,\ldots,\langle x_{n},q\overline{\phi_{1}}(x_{1})\cdots\overline{\phi_{1}}(x_{n-1})\rangle),$$
with $t(x_{n})=s(x_{1})$, which projects to a closed walk $W=(x_{1},\ldots,x_{n})$ in $\Gamma_{1}$ with $\overline{\phi_{1}}(W)=1$.
Since
\begin{align*}
\underline{\phi_{1}}(\overline{W})=\kappa(q)\phi_{1}(W)\kappa(q)^{-1}\in N,
\end{align*}
we have $\phi_{1}(W)=1$ and hence $\phi_{2}(\alpha(W))=1$.
If we write $\overline{\alpha}(u,q)=(\alpha(u),q')$, then
\begin{align*}
\underline{\phi_{2}}(\overline{\alpha}(\overline{W}))=\kappa(q')\phi_{2}(\alpha(W))\kappa(q')^{-1}=1,
\end{align*}
and hence $\overline{\alpha}$ can be lifted to $(\Gamma_{1}\times_{\overline{\phi_{1}}}Q)\times_{\underline{\phi_{1}}}N\to(\Gamma_{2}\times_{\overline{\phi_{2}}}Q)\times_{\underline{\phi_{2}}}N$.

($\Leftarrow$) If $\alpha$ lifts to $\Gamma_{1}\times_{\overline{\phi_{1}}}Q\to\Gamma_{2}\times_{\overline{\phi_{2}}}Q$,
and any such lift lifts again to $(\Gamma_{1}\times_{\overline{\phi_{1}}}Q)\times_{\underline{\phi_{1}}}N
\to(\Gamma_{2}\times_{\overline{\phi_{2}}}Q)\times_{\underline{\phi_{2}}}N$, say $\tilde{\alpha}$, then it can be checked that $\Xi_{\phi_{2}}^{-1}\circ\tilde{\alpha}\circ\Xi_{\phi_{1}}$ is a lift of $\alpha$ to $\Gamma\times_{\phi_{1}}A\to\Gamma\times_{\phi_{2}}A$.

To prove the last assertion, use Proposition~\ref{prop:lifting} again to see that any two lifts of $\alpha$ to $\Gamma_{1}\times_{\overline{\phi_{1}}}Q\to\Gamma_{2}\times_{\overline{\phi_{2}}}Q$, say, $\overline{\alpha}$ and $\overline{\alpha}'$, differ by a left action of $q$ on $\Gamma_{2}\times_{\overline{\phi_{2}}}Q$ for some $q\in Q$. But obviously this action of $q$ as an automorphism can be lifted to the left action of $\kappa(q)$ (actually any $g\in A$ with $\jmath(g)=q$) on $(\Gamma_{2}\times_{\overline{\phi_{2}}}Q)\times_{\underline{\phi_{2}}}N$. Thus whenever $\overline{\alpha}$  lifts again, so does $\overline{\alpha}'$.
\end{proof}

\subsection{Lifting along solvable covers}\label{Sec=SolCov}

Let $\phi:D(\Gamma)\rightarrow A$ be a voltage assignment and let $A$ be a solvable group.
Choose a finite series
\begin{align}
1=A^{(d+1)}\vartriangleleft A^{(d)}\vartriangleleft\cdots\vartriangleleft A^{(0)}=A, \label{eq:chain}
\end{align}
such that, for each $i\in\{0,1,\ldots,d\}$, $A^{(i)}$ is a characteristic subgroup of $A$ and $A^{(i)}/A^{(i+1)}$ is abelian.
Every solvable group $A$ always has such a series, namely, the {\it derived series}, in which $A^{(i+1)}=[A^{(i)},A^{(i)}]$
is the commutator subgroup of $A^{(i)}$ for each $i$.

For each $i\in\{0,1,\ldots,d\}$, choose a section $\kappa^{(i)}:A^{(i)}/A^{(i+1)}\to A^{(i)}$ with $\kappa^{(i)}(1)=1$. For the given voltage assignment $\phi:D(\Gamma)\to A$, we simply use $\Gamma(\phi)$ to denote  $\Gamma\times_{\phi}A$
when the voltage group $A$ is clear from the context.

Set $\Gamma^{(0)}=\Gamma$, $\phi^{(0)}=\phi$. Define an abelian cover $\Gamma^{(i)}\to\Gamma^{(i-1)}$ and a voltage assignment
$\phi^{(i)}:D(\Gamma^{(i)})\rightarrow A^{(i)}$ recursively by setting
\begin{align}
\Gamma^{(i+1)}&=\Gamma^{(i)}(\overline{\phi}^{(i)}), \\
\phi^{(i+1)}\langle x,q\rangle&=\kappa^{(i)}(q)\phi^{(i)}(x)[\kappa^{(i)}
(q\overline{\phi}^{(i)}(x))]^{-1},
\end{align}
where $\overline{\phi}^{(i)}:D(\Gamma^{(i)})\rightarrow A^{(i)}/A^{(i+1)}$ is the composite of $\phi^{(i)}$ with the quotient homomorphism
$A^{(i)}\rightarrow A^{(i)}/A^{(i+1)}$.

For each $i$ there is an isomorphism
\begin{align}
\Xi_{i}:=\Xi_{\phi^{(i)}}:\Gamma^{(i)}(\phi^{(i)})\cong\Gamma^{(i+1)}(\phi^{(i+1)}),
\end{align}
which is an equivalence of $A^{(i)}$-covers of $\Gamma^{(i)}$.

Now we have a commutative diagram
\begin{align}\label{Eq=layer}
\xymatrix{
\Gamma^{(d)}(\phi^{(d)})\ar[d]\ar[r]&\Gamma^{(d-1)}(\phi^{(d-1)})\ar[d]\ar[r]&\ar[r]\cdots\ar[r]&\Gamma^{(1)}(\phi^{(1)})
\ar[d]\ar[r]&\Gamma^{(0)}(\phi^{(0)})\ar[d] \\
\Gamma^{(d)}\ar[r]&\Gamma^{(d-1)}\ar[r]&\ar[r]\cdots\ar[r]&
\Gamma^{(1)}\ar[r]&\Gamma^{(0)}
}
\end{align}
where each small square
$$
\xymatrix{
\Gamma^{(i)}(\phi^{(i)})\ar[d]^{A^{(i)}}\ar[r]^{ \ \ \Xi_{i} \ \ \ \ \ }&\Gamma^{(i-1)}(\phi^{(i-1)})\ar[d]^{A^{(i-1)}} \\
\Gamma^{(i)}\ar[r]^{\ \ \, \, \ \ A^{(i-1)}/A^{(i)} \ \ \, \, \ \ \ \  }&\Gamma^{(i-1)}
},
$$ shows that
the $A^{(i-1)}$-cover $\Gamma^{(i-1)}(\phi^{(i-1)})\to\Gamma^{(i-1)}$ is decomposed into the $A^{(i)}$-cover $\Gamma^{(i)}(\phi^{(i)})\to\Gamma^{(i)}$
and the abelian $(A^{(i-1)}/A^{(i)})$-cover $\Gamma^{(i)}\to\Gamma^{(i-1)}$, and hence the bottom row in the diagram~(\ref{Eq=layer}) is
a layer of abelian covers. Furthermore, in the leftmost square
$$
\xymatrix{
\Gamma^{(d)}(\phi^{(d)})\ar[d]^{A^{(d)}}\ar[r]^{ \ \ \Xi_{d} \ \ \ \ \ }&\Gamma^{(d-1)}(\phi^{(d-1)})\ar[d]^{A^{(d-1)}} \\
\Gamma^{(d)}\ar[r]^{\ \ \, \, \ \ A^{(d-1)}/A^{(d)} \ \ \, \, \ \ \ \  }&\Gamma^{(d-1)}
},
$$
the left vertical projection is also an abelian $A^{(d)}$-cover.
By repeatedly applying Theorem \ref{Thm=lift-decompose}, one can have

\begin{thm} \label{Thm=solvable}
An automorphism $\alpha^{(0)}\in{\rm Aut}(\Gamma)$ can be lifted to a solvable cover $\Gamma\times_{\phi}A$ if and only if
there exist automorphisms
$\alpha^{(i)}\in{\rm Aut}(\Gamma^{(i)})$, $i=1,\ldots,d+1$,
such that, for each $i$, $\alpha^{(i)}$ is a lift of $\alpha^{(i-1)}$ along the abelian cover $\Gamma^{(i)}\to\Gamma^{(i-1)}$.
\end{thm}

\begin{rmk}
\rm Let $G\leqslant\textrm{Aut}(\Gamma)$ be an automorphism group of a graph $\Gamma$. A cover $\tilde{\Gamma}\to\Gamma$ is said to be \emph{$G$-admissible} if each $\alpha\in G$ can be lifted to $\tilde{\Gamma}$. Chen and Shen \cite{CS13} proposed a practical method for finding all $G$-admissible abelian covers and applied it to classify all arc-transitive abelian covers of the Petersen graph. Also, for cubic graphs, Conder and Ma \cite{CM13} develop a method for determining $G$-admissible abelian covers.

Now based on Theorem \ref{Thm=lift-decompose} and Theorem \ref{Thm=solvable}, we are in principle able to find all $G$-admissible $A$-covers for a solvable group $A$ via the technique of step-by-step lifting. This will be illustrated in the next section.
\end{rmk}

\section{Application to metacyclic covers of the tetrahedron branched at face-centers}

\subsection{Preparation}

In this paper, by a {\it map} $\mathcal{M}=(\Gamma,\Sigma)$, we mean a cellular embedding of a graph $\Gamma$ in a closed surface $\Sigma$, where ``cellular" means that, when the graph is deleted, each connected component is homeomorphic to a disk. Call $\Gamma$ and $\Sigma$ the {\it underlying graph} and {\it underlying surface} of $\mathcal{M}$, respectively. There are various other equivalent definitions of maps; see \cite{JS78}.
Given two maps $\mathcal{M}=(\Gamma,\Sigma)$, $\mathcal{M}'=(\Gamma',\Sigma')$,
an isomorphism $\alpha:\mathcal{M}\to\mathcal{M}'$ is a graph isomorphism $\alpha:\Gamma\to\Gamma'$ which can be extended to an orientation-preserving homeomorphism $\Sigma\to\Sigma'$.
The map $\mathcal{M}$ is called {\it regular} if its automorphism group ${\rm Aut}(\mathcal{M})$, which can be identified with a subgroup of
${\rm Aut}(\Gamma)$, acts transitively on $D(\Gamma)$.

In this section, we only consider regular maps. Given maps $\tilde{\mathcal{M}}=(\tilde{\Gamma},\tilde{\Sigma})$, $\mathcal{M}=(\Gamma,\Sigma)$, a {\it regular cover smooth at vertices} $\Pi:\tilde{\mathcal{M}}\to\mathcal{M}$ is a regular cover $\pi:\tilde{\Gamma}\to\Gamma$ which can be extended to a branched cover $\tilde{\Sigma}\to\Sigma$, and $\Pi$ is called an $A$-cover if ${\rm CTG}(\pi)\cong A$. In below, we abbreviate ``regular cover smooth at vertices" to ``regular cover". It is known that regular covers of a map $\mathcal{M}$ are the same as ${\rm Aut}(\mathcal{M})$-admissible regular covers of the underlying graph of $\mathcal{M}$; see Theorem 4.3.3 and Theorem 4.3.5 in \cite{GT87}, and also \cite{MNS02}.

\begin{figure} [h]
  \centering
  \includegraphics[width=0.7\textwidth]{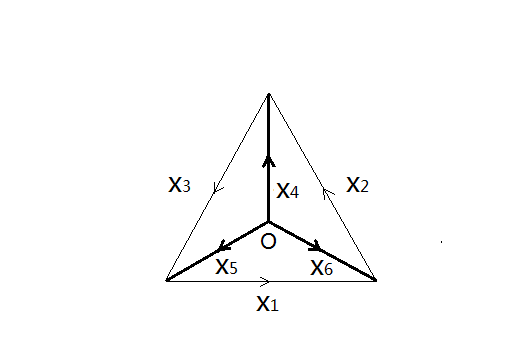}\\
  \caption{The graph $\Gamma$}\label{fig:tetr}
\end{figure}

From now on, let $\Gamma=K_{4}$, as shown in Figure \ref{fig:tetr}. The {\it tetrahedron} $\mathcal{M}(3,3)$ is the standard embedding of $\Gamma$ in the 2-sphere. It is one of the five {\it Platonic maps} (one may refer to \cite{CM92}). For Platonic maps, cyclic covers were classified by Jones and Surowski \cite{JS00}, and elementary abelian covers have been partially classified by Jones \cite{Jo12}.
Now we are going to classify metacyclic covers of $\mathcal{M}(3,3)$. As a related result, 2-arc-transitive metacyclic covers of $K_{4}$ has been classified in Xu et al. \cite{XDKX14}.

It is well-known  (see \cite{MNS02} Example 4.5) that ${\rm Aut}(\mathcal{M}(3,3))$ is isomorphic to the alternating group $\mathbb{A}_{4}$, and can be generated by two automorphisms $\alpha$ and $\beta$ whose actions on arcs are determined as follows:
\begin{align*}
\alpha&:x_{1}\mapsto x_{2}, \ \ x_{2}\mapsto x_{3}, \ \ x_{3}\mapsto x_{1}, \ \ x_{4}\mapsto x_{5}, \ \ x_{5}\mapsto x_{6}, \ \ x_{6}\mapsto x_{4}, \\
\beta&:x_{1}\mapsto x_{1}^{-1}, \ \  x_{2}\mapsto x_{5}^{-1}, \ \ x_{3}\mapsto x_{6},
\ \ x_{4}\mapsto x_{4}^{-1}, \ \ x_{5}\mapsto x_{2}^{-1}, \ \ x_{6}\mapsto x_{3},
\end{align*}
where it is understood that $\alpha(x_{i}^{-1})=(\alpha(x_{i}))^{-1}$, $\beta(x_{i}^{-1})=(\beta(x_{i}))^{-1}$.

A \emph{metacyclic group} is an extension of a cyclic group by another cyclic group. As a well-known fact (see Section 3.7 of Zassenhaus\cite{Za56}), each metacyclic group can be presented as
\begin{align}
H=H(n,m;r,t)=\langle a,b|a^{n}=1, b^{m}=a^{t}, bab^{-1}=a^{r}\rangle. \label{eq:presentation}
\end{align}
for some positive integers $n,m,r,t$ satisfying
\begin{align}
r^{m}-1\equiv t(r-1)\equiv 0\pmod{n}.   \label{eq:condition}
\end{align}
It is easy to see that each element of $H$ can be written as $a^{u}b^{v}$ for some integers $u,v$, and a direct computation shows
\begin{align}
[a^{u_{1}}b^{v_{1}},a^{u_{2}}b^{v_{2}}]&=a^{(r^{v_{1}}-1)u_{2}-(r^{v_{2}}-1)u_{1}}.
\end{align}
Thus the commutator subgroup $[H,H]$ is generated by $a^{d}$, with $d=(r-1,n)$, and the abelianization
$H^{\rm ab}:=H/[H,H]$ is isomorphic to $\mathbb{Z}_{d}\times\mathbb{Z}_{m}$.

For convenience, for any group $A$, we simply call $\langle\alpha,\beta\rangle$-admissible $A$-covers {\it admissible $A$-covers}, and say a voltage assignment is {\it admissible} if its derived cover is admissible.

The task is, for each metacyclic group $H$, to classify admissible $H$-covers of $\Gamma$.
Let $T$ be the spanning tree of $\Gamma$ consisting of the underlying edges of $x_{4},x_{5},x_{6}$. By Theorem \ref{thm:bijection}, Theorem \ref{Thm=lift-decompose} and Remark \ref{rmk:T-reduced}, the procedure can be
divided into two steps:
first determine all $T$-reduced admissible voltage assignments $\varphi:D(\Gamma)\to H^{{\rm ab}}$, next for each such $\varphi$ we determine all voltage assignments $\psi\in\mathcal{C}_{\kappa}(\varphi;[H,H])$, (where $\kappa:H^{{\rm ab}}\to H$ is a section), such that its derived cover admits lifts for certain automorphisms and $\psi\langle x,0\rangle=0$ for all $x\in D(T)$.

\subsection{$H^{{\rm ab}}$-covers}

Clearly $H_{1}(\Gamma;\mathbb{Z})\cong\mathbb{Z}^{3}$ is freely generated by $L_{1}, L_{2}, L_{3}$, with
\begin{align*}
L_{1}=(x_{5},x_{1},x_{6}^{-1}), \ \ \ L_{2}=(x_{6},x_{2},x_{4}^{-1}), \ \ \ L_{3}=(x_{4},x_{3},x_{5}^{-1}).
\end{align*}
The actions induced by $\alpha, \beta$ on $H_{1}(\Gamma;\mathbb{Z})$ are given by
\begin{align*}
\alpha_{\ast}&:L_{1}\mapsto L_{2}, \ \ \ \ L_{2}\mapsto L_{3}, \ \ \ \ L_{3}\mapsto L_{1}, \\
\beta_{\ast}&: L_{1}\mapsto -L_{1}-L_{2}-L_{3}, \ \ \ \ L_{2}\mapsto L_{3}, \ \ \ \ L_{3}\mapsto L_{2},
\end{align*}
hence the matrices representing them, determined via (\ref{eq:matrix}), are given by
\begin{align}
[\alpha]=\left(\begin{array}{ccc} 0 & 1 & 0 \\ 0 & 0 & 1 \\ 1 & 0 & 0 \\ \end{array}\right), \ \ \
[\beta]=\left(\begin{array}{ccc} -1 & -1 & -1 \\ 0 & 0 & 1 \\ 0 & 1 & 0 \\ \end{array}\right).
\end{align}

\begin{lem} \label{lem:p-group}
Let $p$ be a prime number.

{\rm(i)} If $\varphi:D(\Gamma)\to\mathbb{Z}_{p^{k}}$ is a $T$-reduced admissible voltage assignment, then $p^{k}\in\{2,4\}$ and $\varphi(x_{1})=\varphi(x_{2})=\varphi(x_{3})$.

{\rm(ii)} For any $k,k'\geq 1$, there does not exist an admissible voltage assignment $\varphi:D(\Gamma)\to\mathbb{Z}_{p^{k}}\times\mathbb{Z}_{p^{k'}}$.
\end{lem}

\begin{proof}
(i) Let $\varphi:D(\Gamma)\to \mathbb{Z}_{p^{k}}$ be a $T$-reduced admissible voltage assignment. Since $\varphi(x_{1}), \varphi(x_{2}), \varphi(x_{3})$ generate $\mathbb{Z}_{p^{k}}$, at least one of them is invertible; this implies that all of them are invertible, as $\alpha_{\ast}$ cyclically permutes $\varphi(x_{1}), \varphi(x_{2}), \varphi(x_{3})$. Let $u=\varphi(x_{2})\varphi(x_{1})^{-1}$, $v=\varphi(x_{3})\varphi(x_{1})^{-1}$.
To convert $[\varphi]$ into normal form, we take
$B=\left(\begin{array}{ccc} 1 & 0 & 0 \\ -u & 1 & 0 \\ -v & 0 & 1 \\ \end{array}\right)$ and
$C=\left(\begin{array}{ccc} \varphi(x_{1})^{-1} & 0 & 0 \\ 0 & 1 & 0 \\ 0 & 0 & 1 \\ \end{array}\right),$
so that $B[\varphi]C=\left(\begin{array}{ccc} 1 \\ 0 \\ 0 \\ \end{array}\right)$, $s_{1}=0, s_{2}=s_{3}=k$, in the notation of Section 2.

Computing directly, (obscuring the elements we are not interested in)
$$B[\alpha]B^{-1}=\left(\begin{array}{ccc} \ast & \ast & \ast \\ v-u^{2} & \ast & \ast \\ 1-uv & \ast & \ast \\ \end{array}\right), \ \ \
B[\beta]B^{-1}=\left(\begin{array}{ccc} \ast & \ast & \ast \\ (1+u)(u+v) & \ast & \ast \\ (1+v)(u+v) & \ast & \ast \\ \end{array}\right).$$
Applying Theorem \ref{thm:abelian}, we obtain the following equations in $\mathbb{Z}_{p^{k}}$:
$$v-u^{2}=1-uv=0, \ \ \ \ \  (1+u)(u+v)=(1+v)(u+v)=0.$$
They have a solution only when $p^{k}\in\{2,4\}$, and then $u=v=1$.

(ii) First consider the special case $k=k'=1$. Let $\varphi:D(\Gamma)\to \mathbb{Z}_{p}^{2}$ be a $T$-reduced admissible voltage assignment, with $\varphi(x_{i})=(u_{i},v_{i}), 1\leq i\leq 3$.
Since $\varphi(x_{1}), \varphi(x_{2}), \varphi(x_{3})$ generate $\mathbb{Z}_{p}^{2}$, at least one of
$$\left|\begin{array}{cc} u_{1} & v_{1} \\ u_{2} & v_{2} \\ \end{array}\right|, \left|\begin{array}{cc} u_{2} & v_{2} \\ u_{3} & v_{3} \\ \end{array}\right|,
\left|\begin{array}{cc} u_{3} & v_{3} \\ u_{1} & v_{1} \\ \end{array}\right|$$
is invertible. Then all of them are invertible, as $\alpha_{\ast}$ cyclically permutes $\varphi(x_{1}), \varphi(x_{2}), \varphi(x_{3})$.
We may find $B=\left(\begin{array}{ccc} 1 & 0 & 0 \\ u & 1 & 0 \\ v & w & 1 \\ \end{array}\right)$ for some $u,v,w$ and another invertible matrix $C$ such that $B[\varphi]C=\left(\begin{array}{ccc} 1 & 0 \\ 0 & 1 \\ 0 & 0 \\ \end{array}\right)$, and $s_{1}=s_{2}=0, s_{3}=1$.

Computing directly,
\begin{align*}
B[\alpha]B^{-1}&=\left(\begin{array}{ccc} \ast & \ast & \ast \\ \ast & \ast & \ast \\ 1-vw+u(w^{2}-v) & v-w^{2} & \ast \\ \end{array}\right), \\
B[\beta]B^{-1}&=\left(\begin{array}{ccc} \ast & \ast & \ast \\ \ast & \ast & \ast \\ u(v-w-1) & 1+w-v & \ast \\ \end{array}\right).
\end{align*}
By Theorem \ref{thm:abelian} the following hold in $\mathbb{Z}_{p}$:
\begin{align*}
1-vw+u(w^{2}-v)=v-w^{2}&=0, \\
u(v-w-1)=1+w-v&=0.
\end{align*}
But these equations have no solution, as can be checked.

Thus $\varphi:D(\Gamma)\to\mathbb{Z}_{p}^{2}$ cannot be admissible.

Now for general case, let $\varphi:D(\Gamma)\to\mathbb{Z}_{p^{k}}\times\mathbb{Z}_{p^{k'}}$ be a voltage assignment, with $k,k'\geq 1$. The kernel of the quotient homomorphism $\varpi:\mathbb{Z}_{p^{k}}\times\mathbb{Z}_{p^{k'}}\to\mathbb{Z}_{p}^{2}$ is the subgroup generated by $(p,0)$ and $(0,p)$, hence it is characteristic.
If $\varphi$ is admissible, then by Theorem \ref{Thm=lift-decompose}, the composite $\varpi\circ\varphi:D(\Gamma)\to\mathbb{Z}_{p}^{2}$ would also be admissible, contradicting to the above result. Thus there does not exist an admissible voltage assignment
$\varphi:D(\Gamma)\to\mathbb{Z}_{p^{k}}\times\mathbb{Z}_{p^{k'}}$.
\end{proof}
\begin{rmk} \label{rmk:key}
\rm The result (i) recovers part of Proposition 4 of \cite{JS00}, and (ii) generalizes the result obtained in \cite{Jo12} Section 3 that there does not exist an admissible $\mathbb{Z}_{p}^{2}$-cover when $p>2$.

As a consequence, if $\varphi:D(\Gamma)\to A$ is admissible with $A$ an abelian group, then $A\cong\mathbb{Z}_{2}$ or $A\cong\mathbb{Z}_{4}$. This is because,  for any Sylow $p$-subgroup $A_{p}$ of $A$, the kernel of the projection $A\twoheadrightarrow A_{p}$ is characteristic.
\end{rmk}

If there exists an admissible $H$-cover with $H$ a metacyclic group, then $H^{{\rm ab}}\cong\mathbb{Z}_{2^{k}}$ with $k\in\{1,2\}$, hence $2\nmid n$, $(r-1,n)=1$ and
\begin{align}
H\cong H(n,2^{k};r):=\langle a,b|a^{n}=b^{2^{k}}=1,b^{-1}ab=a^{r}\rangle.
\end{align}
For $\eta=\pm 1\in\mathbb{Z}_{2^{k}}$, we may construct an admissible voltage assignment
$$\varphi_{\eta}:D(\Gamma)\to\mathbb{Z}_{2^{k}}, \ \ \ \ \ x_{i}\mapsto \eta, \ \ \ x_{i+3}\mapsto 0, \ \ \ (1\leq i\leq 3).$$
Lift $\alpha,\beta$ to $\tilde{\alpha},\tilde{\beta}\in{\rm Aut}(\Gamma\times_{\varphi_{\eta}}\mathbb{Z}_{2^{k}})$, respectively, with
\begin{align}
\tilde{\alpha}\langle x_{i},g\rangle=\langle \alpha(x_{i}),g\rangle,  \ \ \ \ \
\tilde{\beta}\langle x_{i},g\rangle=\langle \beta(x_{i}),g\rangle, \ \ \ \ (1\leq i\leq 6).
\end{align}

\subsection{$[H,H]$-covers of $\Gamma\times_{\varphi}H^{\rm ab}$}

Denote $m=2^{k}, k\in\{1,2\}$, and for $\eta\in\{\pm 1\}$ abbreviate $\Gamma\times_{\varphi_{\eta}}\mathbb{Z}_{m}$ to $\Gamma(\varphi_{\eta})$.

Choose a section
$$\kappa:\mathbb{Z}_{m}\to H(n,m;r), \ \ \ \ \ g\mapsto b^{g}$$
as a right inverse to the quotient homomorphism $H(n,m;r)\twoheadrightarrow\mathbb{Z}_{m}$.

With a characteristic subgroup $\mathbb{Z}_{n}$, we shall find voltage assignments
$\psi\in\mathcal{C}_{\kappa}(\varphi_{\eta};\mathbb{Z}_{n})$ such that $\tilde{\alpha}$ and $\tilde{\beta}$ can be lifted to $\Gamma(\varphi_{\eta})\times_{\psi}\mathbb{Z}_{n}$, and
\begin{align}
\langle\psi_{1},\psi_{2},\psi_{3}\rangle=\mathbb{Z}_{n}, \ \ \ \ \ \psi_{4}=\psi_{5}=\psi_{6}=0, \label{eq:psi}
\end{align}
where $\psi_{i}=\psi\langle x_{i},0\rangle$.
The condition (\ref{Eq=well psi}) turns out to be
$$\psi\langle x_{i},g\rangle=\kappa(g)\psi_{i}\kappa(g)^{-1}=r^{g}\psi_{i}, \ \  (1\leq i\leq 6),$$
hence $\psi\langle x_{i},g\rangle=0$ for $i=4,5,6$.

\begin{lem}   \label{lem:lift2}
The automorphisms $\tilde{\alpha}$ and $\tilde{\beta}$ can be lifted if and only if
\begin{align}
u(\psi_{3}-\psi_{2})-v(\psi_{2}-\psi_{1})&=0,  \label{eq:uv1} \\
u(\psi_{1}+\psi_{2}+(1+r^{\eta})\psi_{3})+v((\psi_{1}+\psi_{2})+r^{\eta}(\psi_{2}+\psi_{3}))&=0, \label{eq:uv2}
\end{align}
for any integers $u,v$ satisfying
\begin{align}
u(\psi_{2}-\psi_{1})+v(\psi_{3}-\psi_{1})=0. \label{eq:coeff2}
\end{align}
\end{lem}
\begin{proof}
Let 
$$\tilde{L}_{0}=(\langle x_{1},0\rangle, \langle x_{6},0\rangle^{-1}, \langle x_{5},0\rangle,\ldots,
\langle x_{1},m-1\rangle,\langle x_{6},m-1\rangle^{-1},\langle x_{5},m-1\rangle),$$
and for each $g\in\mathbb{Z}_{m}$, let
\begin{align*}
\tilde{L}_{2,g}&=(\langle x_{2},g\rangle,\langle x_{4},g+\eta\rangle^{-1},\langle x_{6},g+\eta\rangle,\langle x_{1},g\rangle^{-1},\langle x_{5},g\rangle^{-1},\langle x_{6},g\rangle), \\
\tilde{L}_{3,g}&=(\langle x_{3},g\rangle,\langle x_{5},g+\eta\rangle^{-1},\langle x_{6},g+\eta\rangle,\langle x_{1},g\rangle^{-1},\langle x_{5},g\rangle^{-1},\langle x_{4},g\rangle).
\end{align*}
The closed walks $\tilde{L}_{0}, \tilde{L}_{2,g}, \tilde{L}_{3,g}, g\in\mathbb{Z}_{m}$ form a basis of $H_{1}(\Gamma(\varphi_{\eta});\mathbb{Z})\cong\mathbb{Z}^{2m+1}$.
To see this, let $\tilde{T}$ denote the spanning tree of $\Gamma(\varphi_{\eta})$ consisting of the underlying edges of $\langle x_{i},0\rangle$,
$4\leq i\leq 6, g\in\mathbb{Z}_{m}$ and $\langle x_{1},g\rangle$, $g\in\mathbb{Z}_{m}-\{0\}$, then regarded as elements of $H_{1}(\Gamma(\varphi_{\eta});\mathbb{Z})$, each fundamental closed walk can be written as a linear combination of $\tilde{L}_{0}, \tilde{L}_{2,g}, \tilde{L}_{3,g}, g\in\mathbb{Z}_{m}$.

Computing directly, we get
\begin{align*}
\psi(\tilde{L}_{0})=0, \ \ \ \ \ \ \ \psi(\tilde{L}_{2,g})&=r^{g}(\psi_{2}-\psi_{1}), \ \ \ \ \ \ \ \psi(\tilde{L}_{3,g})=r^{g}(\psi_{3}-\psi_{1}). \\
\psi(\tilde{\alpha}(\tilde{L}_{0}))=0, \ \ \   \psi(\tilde{\alpha}(\tilde{L}_{2,g}))&=r^{g}(\psi_{3}-\psi_{2}),  \ \ \ \psi(\tilde{\alpha}(\tilde{L}_{3,g}))=-r^{g}(\psi_{2}-\psi_{1}), \\
\psi(\tilde{\beta}(\tilde{L}_{0}))=0, \ \ \   \psi(\tilde{\beta}(\tilde{L}_{2,g}))&=r^{g}(\psi_{3}+\psi_{2}+\psi_{1})+r^{g+\eta}\psi_{3}, \\ \psi(\tilde{\beta}(\tilde{L}_{3,g}))&=r^{g}(\psi_{2}+\psi_{1})+r^{g+\eta}(\psi_{2}+\psi_{3}).
\end{align*}

For a general element $\tilde{L}=c\tilde{L}_{0}+\sum\limits_{g=1}^{m}(u_{g}\tilde{L}_{2,g}+v_{g}\tilde{L}_{3,g})\in H_{1}(\Gamma(\varphi_{\eta});\mathbb{Z})$,
\begin{align*}
\psi(\tilde{L})=&\left(\sum\limits_{g=1}^{m}u_{g}r^{g}\right)(\psi_{2}-\psi_{1})+\left(\sum\limits_{g=1}^{m}v_{g}r^{g}\right)(\psi_{3}-\psi_{1}),  \\
\psi(\tilde{\alpha}(\tilde{L}))=&\left(\sum\limits_{g=1}^{m}u_{g}r^{g}\right)(\psi_{3}-\psi_{2})-
\left(\sum\limits_{g=1}^{m}v_{g}r^{g}\right)(\psi_{2}-\psi_{1}), \\
\psi(\tilde{\beta}(\tilde{L}))=&\left(\sum\limits_{g=1}^{m}u_{g}r^{g}\right)(\psi_{1}+\psi_{2}+(1+r^{\eta})\psi_{3})+ \\
&\left(\sum\limits_{g=1}^{m}v_{g}r^{g}\right)((\psi_{1}+\psi_{2})+r^{\eta}(\psi_{2}+\psi_{3})).
\end{align*}
A sufficient and necessary condition for $\tilde{\alpha}$ to be lifted to $\Gamma(\varphi_{\eta})\times_{\psi}\mathbb{Z}_{n}$ is that $\psi(\tilde{\alpha}(\tilde{L}))=0$ whenever $\psi(\tilde{L})=0$, and similarly for $\tilde{\beta}$.
Noting that $\sum\limits_{g=1}^{m}u_{g}r^{g}$ and $\sum\limits_{g=1}^{m}v_{g}r^{g}$ may take arbitrary values, we establish the lemma.
\end{proof}

Taking integers $u,v$ with $u\equiv\psi_{3}-\psi_{1}, v\equiv\psi_{1}-\psi_{2}\pmod{n}$ so that (\ref{eq:coeff2}) holds, we deduce from (\ref{eq:uv1}) and (\ref{eq:uv2}) the following identities in the ring $\mathbb{Z}_{n}$:
\begin{align}
(\psi_{2}-\psi_{1})^{2}&=(\psi_{1}-\psi_{3})(\psi_{3}-\psi_{2}),  \label{eq:C1} \\
(r^{\eta}-1)\psi_{2}(\psi_{3}-\psi_{1})&=(r^{\eta}+1)(\psi_{3}^{2}-\psi_{2}^{2}). \label{eq:C2}
\end{align}

Suppose $n=\prod\limits_{p\in\Lambda}p^{k_{p}}$ with $2\notin\Lambda$. Fix $p\in\Lambda$.

If $\psi_{1}\equiv\psi_{2}\equiv\psi_{3}\pmod{p}$, then taking $u=2n/p$ and $v=-n/p$ which satisfy (\ref{eq:coeff2}), one would deduce from (\ref{eq:uv2}) that $4\psi_{1}\equiv 0\pmod{p}$, but this contradicts to (\ref{eq:psi}). Hence
\begin{align}
\psi_{2}\not\equiv\psi_{1} \ \ \ \ \text{or\ \ \ } \psi_{3}\not\equiv\psi_{2}\pmod{p}. \label{eq:coprime}
\end{align}
In the following discussions, by writing $z=z'$, we mean $z\equiv z'\pmod{p^{k_{p}}}$. The conditions $r^{m}=1$ and $r\not\equiv 1\pmod{p}$ imply either $r=-1$ or $r^{2}=-1$.
\begin{enumerate}
  \item[\rm(i)] If $r=-1$, then (\ref{eq:C2}) implies $\psi_{2}(\psi_{3}-\psi_{1})=0$. By (\ref{eq:C1}) and (\ref{eq:coprime}), it is impossible
       that $\psi_{3}-\psi_{1}\equiv 0\pmod{p}$, so $\psi_{2}=0$.
       By (\ref{eq:C1}),
       $$(\psi_{3}-2\psi_{1})^{2}+3\psi_{3}^{2}=0.$$

       If $p=3$, then $k_{3}=1$, because if $k_{3}>1$, then $\psi_{3}-2\psi_{1}\equiv 0\pmod{3}$, which would imply $\psi_{1}\equiv 0\pmod{3}$ and also $\psi_{3}\equiv 0\pmod{3}$, contradicting to (\ref{eq:psi}). Hence $\psi_{3}=2\psi_{1}$.

       If $p>3$, then the Legendre symbol $(-3/p)$ is equal to $1$. By Law of Quadratic Reciprocity (see \cite{IR90} P.53, Theorem 1),
       $$\left(\frac{p}{3}\right)=(-1)^{\frac{p-1}{2}}\cdot\left(\frac{3}{p}\right)=\left(\frac{-1}{p}\right)\left(\frac{3}{p}\right)
       =\left(\frac{-3}{p}\right)=1,$$
       which is equivalent to $p\equiv 1\pmod{3}$. On the other hand, by \cite{IR90} Proposition 4.2.3, this is also sufficient for the existence of a square root of $-3$ in $\mathbb{Z}_{p^{k_{p}}}$. Thus $2\psi_{3}=(1+\theta_{p})\psi_{1}$, with $\theta_{p}$ a square root of $-3$.
  \item[\rm(ii)] If $r^{2}=-1$, then $m=4$ and $r^{\eta}=r\eta$. 
       By (\ref{eq:C2}),
       \begin{align}
       \psi_{1}\psi_{2}=\psi_{3}\psi_{2}+r\eta(\psi_{3}^{2}-\psi_{2}^{2}),   \label{eq:psi123}
       \end{align}
       and consequently,
       \begin{align*}
       \psi_{2}^{2}(\psi_{2}-\psi_{1})^{2}=(\psi_{1}\psi_{2}-\psi_{2}^{2})^{2}&=(\psi_{3}-\psi_{2})^{2}(\psi_{2}+r\eta(\psi_{3}+\psi_{2}))^{2},  \\
       \psi_{2}^{2}(\psi_{3}-\psi_{2})(\psi_{1}-\psi_{3})&=r\eta(\psi_{3}-\psi_{2})^{2}\psi_{2}(\psi_{3}+\psi_{2}).
       \end{align*}
       By (\ref{eq:C1}), the right-hand-sides of these two equations coincide, 
       and by (\ref{eq:C1}), (\ref{eq:coprime}), $\psi_{3}-\psi_{2}\not\equiv 0\pmod{p}$, hence
       \begin{align}
       z^{2}+\psi_{2}z+\psi_{2}^{2}=0, \ \ \ \ \ \text{with\ } z=r\eta(\psi_{3}+\psi_{2}). \label{eq:last}
       \end{align}
       By (\ref{eq:psi}), at least one of $\psi_{1},\psi_{2},\psi_{3}$ is invertible, this together with (\ref{eq:C1}) and (\ref{eq:last}) implies that $\psi_{2}$ is invertible. Thus (\ref{eq:last}) has a solution if and only if $p=3$ or $p\equiv 1\pmod{3}$, but $p=3$ is impossible
       since $r^{2}=-1$. Now
       $2z=(\theta_{p}-1)\psi_{2}$, with $\theta_{p}$ a square root of $-3$ in $\mathbb{Z}_{p^{k_{p}}}$,
       hence by (\ref{eq:last}), $2\psi_{3}=(r\eta(1-\theta_{p})-2)\psi_{2}$,
       and then by (\ref{eq:psi123}), $\psi_{1}=(r\eta-\theta_{p})\psi_{2}$.
       Rewrite $\psi_{2},\psi_{3}$ in terms of $\psi_{1}$:
       $$\psi_{2}=(r\eta-\theta_{p})^{-1}\psi_{1}, \ \ \ \ \ \psi_{3}=(1-\theta_{p})^{-1}(r\eta-1)\psi_{1}.$$
\end{enumerate}

We have found all voltage assignments $\psi:D(\Gamma(\varphi_{\eta}))\to\mathbb{Z}_{n}$ satisfying (\ref{eq:psi}) and the conditions in Lemma \ref{lem:lift2}. By (\ref{eq=def(phi)}), the corresponding voltage assignment $\phi:D(\Gamma)\to H$ is given by
$$\phi(x_{i})=a^{\psi_{i}}b^{\eta}, \ \ (1\leq i\leq 3), \ \ \ \ \phi(x_{4})=\phi(x_{5})=\phi(x_{6})=1.$$
Note that $\psi_{1}$ is invertible in $\mathbb{Z}_{n}$, and there exists an automorphism of $H$ taking $a$ to $a^{\psi_{1}^{-1}}$ and fixing $b$,
so we may assume $\psi_{1}=1$.

\subsection{The result}

Given $\eta\in\{\pm 1\}$, and $\theta=\{\theta_{p}\}_{p\in\Lambda}$ with $\theta_{p}$ a square root of $-3$ in $\mathbb{Z}_{p^{k_{p}}}$,
($\theta_{p}=0$ if $p=3$).
For each $p\in\Lambda$, put
\begin{align}
\mu(\eta,\theta_{p})&=
\left\{
\begin{array}{cc}
0, \ \ \ \ &r\equiv -1\pmod{p^{k_{p}}}, \\
(r\eta-\theta_{p})^{-1}, \ \ &r^{2}\equiv -1\pmod{p^{k_{p}}},
\end{array}\right.  \\
\nu(\eta,\theta_{p})&=
\left\{
\begin{array}{cc}
(1+\theta_{p})/2, \ \ \ &r\equiv -1\pmod{p^{k_{p}}}, \\
(1-\theta_{p})^{-1}(r\eta-1), \  \ &r^{2}\equiv -1\pmod{p^{k_{p}}}.
\end{array}\right.
\end{align}
Let $\mu(\eta,\theta)$ and $\nu(\eta,\theta)$ denote the image of $\{\mu(\eta,\theta_{p})\}_{p\in\Lambda}$ and $\{\nu(\eta,\theta_{p})\}_{p\in\Lambda}$, respectively, under the canonical isomorphism $\prod\limits_{p\in\Lambda}\mathbb{Z}_{p^{k_{p}}}\cong\mathbb{Z}_{n}$,
and define the voltage assignment $\phi^{\eta,\theta}:D(\Gamma)\to H$ by
\begin{align}
x_{1}\mapsto ab^{\eta}, \ \ x_{2}\mapsto a^{\mu(\eta,\theta)}b^{\eta}, \ \ x_{3}\mapsto a^{\nu(\eta,\theta)}b^{\eta}, \ \ \
x_{i}\mapsto 1, \  (4\leq i\leq 6).
\end{align}

\begin{lem}
There exists an automorphism of $H$ taking $ab$ to $ab^{-1}$ if and only if $r\equiv -1\pmod{n}$.
\end{lem}
\begin{proof}
If $r\equiv -1\pmod{n}$, then it is easy to verify that $\sigma:H\to H$, $a^{u}b^{v}\mapsto a^{u}b^{-v}$ is a well-defined automorphism with $\sigma(ab)=ab^{-1}$.

Conversely, suppose $\sigma$ is an automorphism of $H$ with $\sigma(ab)=ab^{-1}$. Since $\langle a\rangle$ is characteristic, we have $\sigma(a)=a^{u}$ for some $u$ coprime to $n$, hence $\sigma(b)=a^{1-u}b^{-1}$. It follows from
$$(a^{1-u}b^{-1})a^{u}(ba^{u-1})=\sigma(b)\sigma(a)\sigma(b)^{-1}=\sigma(bab^{-1})=\sigma(a^{r})=a^{ru}$$
that $a^{u}=ba^{ru}b^{-1}=(bab^{-1})^{ru}=a^{r^{2}u}$. Thus $r\equiv -1\pmod{n}$, (recalling that $(r-1,n)=1$).
\end{proof}

Therefore, by Proposition \ref{prop:lifting}, when $r\equiv -1\pmod{n}$, $\Gamma\times_{\phi^{-1,\theta}}H$ is equivalent to $\Gamma\times_{\phi^{1,\theta}}H$, but when $r\not\equiv -1\pmod{n}$, $\Gamma\times_{\phi^{\eta,\theta}}H$ is inequivalent to each other for different pairs $(\eta,\theta)$.

\medskip

\begin{nota}
\rm Let $\Omega$ denote the set of odd natural numbers $n$ such that $9\nmid n$ and $p\equiv 1\pmod{3}$ for each prime divisor $p$ of $n$ with $p>3$.

For $n\in\Omega$, let $o(n)$ denote the number of its prime divisors larger than $3$.
\end{nota}

\begin{thm}
{\rm(i)} A metacyclic group $H$ is the covering transformation group of some regular cover of $\mathcal{M}(3,3)$ if and only if
\begin{itemize}
  \item $H\cong H(n,2^{k};-1)$ with $k\in\{1,2\}$ and $n\in\Omega$; or
  \item $H\cong H(n,4;r)$ with $n\in\Omega$, $(r-1,n)=1$, $n\nmid r+1$ and $n\mid r^{4}-1$.
\end{itemize}

{\rm(ii)} Each $H(n,2^{k};-1)$-cover of $\mathcal{M}(3,3)$ is equivalent to $\Gamma\times_{\phi^{1,\theta}}H(n,2^{k};-1)$ for a unique $\theta$.
There are $2^{o(n)}$ $H(n,2^{k};-1)$-covers in total.

Each $H(n,4;r)$-cover of $\mathcal{M}(3,3)$ is equivalent to $\Gamma\times_{\phi^{\eta,\theta}}H(n,2^{k};r)$ for a unique pair $(\eta,\theta)$.
There are $2^{o(n)+1}$ $H(n,4;r)$-covers in total.
\end{thm}

\vspace{1cm}

{\bf Acknowledgements:}

J.H. Kwak is partially supported by KRF-2012007478.



\begin{thebibliography}{}


\bibitem{Ch13}
H.M. Chen,
\textsl{Lifting automorphisms along abelian regular coverings of graphs}.
Discrete Math. 313 (2013),  1535--1539.


\bibitem{CS13}
H.M. Chen, H. Shen,
\textsl{How to find $G$-admissible coverings of a graph?}.
Linear Algebra Appl. 438 (2013), 3303--3320.


\bibitem{CM13}
M.D.E. Conder, J.C. Ma,
\textsl{Arc-transitive abelian regular covers of cubic graphs}.
J. Algebra 387 (2013),  215--242.


\bibitem{CM92}
R. Cori, A. Mach\`{i},
\textsl{Maps, hypermaps and their automorphisms: a survey, I, II, III}.
Expositiones Math. 10 (1992), 403--427, 429--447, 449--467.



\bibitem{DKX03}
S.F. Du, J.H. Kwak, M.Y. Xu,
\textsl{Linear criteria for lifting automorphisms of abelian regular coverings}.
Linear Algebra Appl. 373 (2003), 101--119.


\bibitem{GT77}
J.L. Gross, T.W. Tucker,
\textsl{Generating all graph coverings by permutation voltage assignments}.
Discrete Math. 18 (1977), 273--283.


\bibitem{GT87}
J.L. Gross, T.W. Tucker,
\textsl{Topological graph theory}.
A Wiley-Interscience Publication, 1987.



\bibitem{HKL96}
S. Hong,  J.H. Kwak,  J. Lee,
\textsl{Regular graph coverings whose covering transformation groups have the isomorphism extension property}.
Discrete Math. 148 (1996), 85--105.


\bibitem{IR90}
K. Ireland, M. Rosen,
\textsl{A Classical Introduction to Modern Number Theory}.
Graduate Texts in Mathematics vol. 84, 2nd edition, Springer-Verlag, New York, 1990.



\bibitem{Jo12}
G.A. Jones,
\textsl{Elementary abelian regular coverings of Platonic maps, Case I: ordinary representations}.
J. Algebra Combin. 41 (2015), no. 2, 461--491.


\bibitem{JS78}
G.A. Jones, D. Singerman,
\textsl{Theorey of maps on orientable surfaces}.
Proc. London Math. Soc. 37 (1978), 273--307.



\bibitem{JS00}
G.A. Jones, D.B. Surowski,
\textsl{Regular Cyclic Coverings of the Platonic Maps}.
European. J. Combinatorics 21 (2000), 333--345.


\bibitem{Ma98}
A. Malni\v c,
\textsl{Group actions, coverings and lifts of automorphisms}.
Discrete Math. 182 (1998), 203--218.

\bibitem{MMP-cover2004}
A. Malni\v c, Maru\v si\v c, P. Poto\v cnik,
\textsl{Elementary abelian covers of graphs},
J. Algebraic Combin. 20 (2004), 71--97.


\bibitem{MMP2004}
A. Malni\v c, D. Maru\v si\v c, P. Poto\v cnik,
\textsl{On cubic graphs admitting an edge-transitive solvable group},
J. Algebraic Combin. 20 (2004), 99--113.

\bibitem{MNS00}
A. Malni\v c, R. Nedela, M. \v Skoviera,
\textsl{Lifting graph automorphisms by voltage assignments.}
European J. Combin. 21 (2000), 927--947.


\bibitem{MNS02}
A. Malni\v c, R. Nedela, M. \v Skoviera,
\textsl{Regular homomorphisms and regular maps.}
European J. Combin. 23 (2002), 449--461.


\bibitem{MP06}
A. Malni\v c, P. Poto\v cnik,
\textsl{Invariant subspaces, duality, and covers of the Petersen graph}.
European J. Combin 27 (2006), 971--989.


\bibitem{Sk86}
M. \v Skoviera,
\textsl{A contribution to the theory of voltage graphs.}
Discrete Math. 61 (1986), 281--292.


\bibitem{XDKX14}
W.Q. Xu, S.F. Du, J.H. Kwak, M.Y. Xu,
\textsl{2-Arc-transitive metacyclic covers of complete graphs}.
J. Combin. Theory Ser. B 111 (2015), 54--74.


\bibitem{Za56}
H.J. Zassenhaus,
\textsl{The theory of groups}.
Second edition, Chelsea, New York, 1956.


\end{thebibliography}
\end{document}